# Investigations of Mertens and Liouville summation functions

VICTOR VOLFSON

ABSTRACT. Summation arithmetic functions of Mertens and Liouville are investigated in the paper. It is proved that the limiting distribution of these functions is the normal. It is also shown that the estimating of standard deviation of these functions $O(n^{1/2})$ cannot be improved. The estimate of the average of Liouville summation function is found. The estimate of the order of growth of the ratio of the summation functions of Mertens and Liouville is also found.

1. INTRODUCTIN

An arithmetic function in the general case is a function defined on the set of natural numbers and taking values on the set of complex numbers. The name arithmetic function is related to the fact that this function expresses some arithmetic property of the natural series. Many arithmetic functions, considered in number theory, take on natural values. In this paper we will consider just such arithmetic functions.

An arithmetic function of a summation type is a function:

$$S(x) = \sum_{n \leq x} f(n).$$

Examples of summation functions are Mertens summation function $M(x) = \sum_{k \leq x} \mu(k)$, where $\mu(k)$ is Möbius function, and the summation function of Liouville $L(x) = \sum_{k \leq x} \lambda(k)$, where $\lambda(k)$ is Liouville function.

Liouville function $\lambda(k) = 1$ if the natural number k has an even number of prime divisors of the first degree with allowance for the degree and $\lambda(k) = -1$ if the natural number k has an odd number of prime divisors with allowance for the degree.

---

Mobius function $\mu(k)=1$ if the natural number k has an even number of prime divisors of the first degree, $\mu(k)=-1$ if the natural number k has an odd number of prime divisors of the first degree and $\mu(k)=0$ if the natural number k has prime divisors of not only first degree.

The density of prime numbers decreases with increasing n, therefore at the beginning of the natural series more numbers have one prime divisor and, correspondingly, an odd number of prime divisors. Therefore, functions $\lambda(k)$ and $\mu(k)$ at the beginning of the natural series often take negative values and the summation functions of Mertens $M(x)=\sum_{k\leq x}\mu(k)$ and Liouville $L(x)=\sum_{k\leq x}\lambda(k)$ are shifted to the negative range.

Moreover, the summation function of Liouville decreases faster than Mertens summation function, since Möbius function $\mu(k)$ takes on zero values where Liouville function $\lambda(k)$ takes negative values (when $\mu(k)=0$ the natural number k has an odd number of prime divisors of degree higher than the first).

Therefore, it is of interest to study the above summation arithmetic functions and determine the similarity and difference in their behavior.

We will explore in the paper:

1. The limit distribution of the data of the summation arithmetic functions.
2. The asymptotic behavior of the above summation arithmetic functions.

We will use both probabilistic and exact methods for these purposes.

2. LIMIT DISTRIBUTION OF MERTENS AND LIOUVILLE SUMMATION ARITHMETIC FUNCTIONS

It is known [1] that any initial segment of a natural number $\{1,2,...,n\}$ can be transformed into a probability space $(\Omega_n, \mathcal{A}_n, \mathbb{P}_n)$ in a natural way, taking as $\Omega_n=\{1,2,...,n\}$, $\mathcal{A}_n$ - all subsets and $\mathbb{P}_n(A)=\frac{|A|}{n}$. Then an arbitrary (real) function of a natural argument $f(k)$ (more precisely, its restriction to $\Omega_n$) can be regarded as a random variable $\xi_n$ on this probability space: $\xi_n(k)=f(k), 1\leq k\leq n$. In particular, we can talk about mathematical expectation (average value) $-M[\xi_n]=\dfrac{1}{n}\sum_{k=1}^{n}f(k)$, variance - $D[\xi_n]=\dfrac{1}{n}\sum_{k=1}^{n}|f(k)|^2-|\dfrac{1}{n}\sum_{k=1}^{n}f(k)|^2$, distribution function - $F_{\xi_n}(x)=\dfrac{1}{n}\{k\leqslant n:f(k)\leqslant x\}$ and characteristic function - $\varphi_{\xi_n}(t)=\dfrac{1}{n}\sum_{k=1}^{n}e^{itf(k)}$.

The limit distribution function for Mobius arithmetic function was found in [2]:

$$\lim_{n\to\infty} P(\mu(n)<y)=F(y), \qquad (2.1)$$

where $P(\mu(n)<y)$ is the distribution function of Mobius function - $\mu(n)$, and $F(y)$ is the limiting distribution function for $\mu(n)$, which is equal to:

$$F(y)=\begin{cases} 0, & y<-1 \ ; \\ p, & -1\leq y<0 \ ; \\ 1-p, & 0\leq y<1 \ ; \\ 1, & y\geq 1 \ . \end{cases} \qquad (2.2)$$

where in (2.2) $p=3/\pi^2$.

We denote Mobius function with the limiting distribution - $\mu$.

Corollary 1. Having in mind (2.1) and (2.2), the average value of the limiting distribution of Mobius function is equal to $M[\mu]=0$ and the variance of the limiting distribution of Mobius function is equal to $D[\mu]=2p=6/\pi^2$.

Let us find the limiting distribution for Liouville arithmetic function.

Theorem 1

$$\lim_{n\to\infty} P(\lambda(n)<y)=G(y), \qquad (2.3)$$

where $P(\lambda(n)<y)$ is the distribution function of Liouville arithmetic function - $\lambda(n)$ and $G(y)$ is the limiting distribution function - $\lambda(n)$, which is equal to:

$$G(y)=\begin{cases} 0, & y<-1 \ ; \\ 1/2, & -1\leq y<1 \ ; \\ 1, & y\geq 1 \ . \end{cases} \qquad (2.4)$$

.

Proof

We introduce probability spaces - $\{\Omega_n, A_n, P_n\}$, where $\Omega_n = (1, 2, ..., n)$, $A_n$ is the collection of all subsets of $\Omega_n$ and $P_n = \{v_1(n), v_2(n)\}$, where $v_1(n) + v_2(n) = 1$ and

$$P(\lambda(i) = 1) = v_1(n) = \frac{1}{n}\{i \in \{1, 2, ..., n\} : \lambda(i) = 1\},$$

$$P(\lambda(i) = -1) = v_2(n) = \frac{1}{n}\{i \in \{1, 2, ..., n\} : \lambda(i) = -1\}.$$

We introduce a random variable - $x_n : x_n(i) = \lambda(i), 1 \leq i \leq n$, where the probabilities are: $P(\lambda(i) = 1) = v_1(n), P(\lambda(i) = -1) = v_2(n)$.

We denote the distribution function of the random variable $x_n(i)$ - $G_n(y) = P(x_n(i) < y)$. Then $x_n(i)$ has the distribution function:

$$G_n(y) = \begin{cases} 0, y < -1; \\ v_2(n), -1 \leq y < 1; \\ 1, y \geq 1. \end{cases} \tag{2.5}$$

Based on (2.5) and Remarks 4 on p. 123 [3], the distribution functions $G_n(y)$ converge to the distribution function $G(y)$ at a value $n \to \infty$, as discrete distributions having jumps at the same points.

Consequently, $\lim_{n \to \infty} P(\lambda(n) < y) = G(y)$ that corresponds to (2.4).

We denote Liouville function with the limiting distribution - $\lambda$.

Corollary 2. Having in mind (2.3) and (2.4), the average value of the limiting distribution of Liouville function is equal to $M[\lambda] = 0$ and the variance of the limiting distribution of Liouville function is equal to $D[\lambda] = 1$.

The following upper bound was proved for Mobius and Liouville arithmetic functions in [4]:

$$M[f(i)f(j)] - M[f(i)]M[f(j)] = O(1/n), \quad (2.6)$$

where $M[...]$ is the average value of the arithmetic function enclosed in parentheses, and $f(k)$ is Mobius or Liouville arithmetic function.

The asymptotic independence of these arithmetic functions follows from the estimate (2.6):

$$\lim_{n \to \infty} M[f(i)f(j+n)] = M[f(i)]M[f(j)]. \quad (2.7)$$

Theorem 2

The limiting distribution for Mertens and Liouville summation arithmetic functions $S(n)$ is the normal distribution.

Proof

We denote the characteristic function for the arithmetic function $f(k)$:

$$\varphi_{f(k)}(t) = M[e^{itf(k)}]. \quad (2.8)$$

Let us consider the summation arithmetic function $S(n) = \sum_{k=1}^{n} f(k)$, then:

$$S(n) - M[S(n)] = \sum_{k=1}^{n} (f(k) - M[f(k)]). \quad (2.9)$$

Taking into account (2.9) the following equality holds:

$$\varphi_{\frac{S(n) - M[S(n)]}{\sqrt{D[S(n)]}}}(t) = M[e^{it\frac{f(1) - M[f(1)]}{\sqrt{D[S(n)]}}} \ldots e^{it\frac{f(n) - M[f(n)]}{\sqrt{D[S(n)]}}}]. \quad (2.10)$$

Having in mind (2.7) and (2.10), for any fixed t, we obtain:

$$\varphi_{\frac{S(n) - M[S(n)]}{\sqrt{D[S(n)]}}}(t) = M[e^{it\frac{f(1) - M[f(1)]}{\sqrt{D[S(n)]}}}] \ldots M[e^{it\frac{f(n) - M[f(n)]}{\sqrt{D[S(n)]}}}] = \prod_{n=1}^{m} \varphi_{\frac{f(k) - M[f(k)]}{\sqrt{D[S(n)]}}}(t). \quad (2.11)$$

Based on [4] we obtain:

$$\sqrt{D[S(n)]} = O(n^{1/2}). \quad (2.12)$$

We denote Mobius or Liouville function with the limiting distribution - $f$.

Then, having in mind (2.11) and (2.12) we obtain for any fixed t and $n \to \infty$:

$$\varphi_{\frac{S(n)-M[S(n)]}{\sqrt{D[S(n)]}}}(t) = (\varphi_{\frac{f}{O(n^{1/2})}}(t))^n. \tag{2.13}$$

Based on Theorem 1 on page 296 [5] having in mind (2.13) we obtain with the value $t \to 0$:

$$(\varphi_{\frac{f}{\sqrt{n}}}(t))^n = 1 - \frac{t^2}{2} + o(t^2). \tag{2.14}$$

Therefore, based on (2.14) we obtain for any fixed t and $n \to \infty$:

$$\varphi_{\frac{S(n)-M[S(n)]}{\sqrt{D[S(n)]}}}(t) = [1 - \frac{t^2}{2n} + o(\frac{1}{n})]^n \to e^{-t^2/2}. \tag{2.15}$$

The function $e^{-t^2/2}$ in (2.15) is the characteristic function of the normal distribution with zero mean and unit variance, which proves the required assertion, by the continuity theorem on page 343 [5].

Let us consider the characteristics of this limiting distribution.

It was proved in [4] that the standard deviation of Mertens and Liouville summation functions $S(n)$ has the order:

$$\sigma[S(n)] = O(n^{1/2}). \tag{2.16}$$

We show that the estimate (2.16) cannot be improved.

Suppose that the Riemann hypothesis is true:

$$S(n) = O(n^{1/2} h(n)), \tag{2.17}$$

where $h(n)$ is a slowly growing function.

Then we obtain the estimate in Lemma 1 [4]:

$$\sum_{i=1}^{n} \sum_{j=1}^{n} f(i)f(j) = o(nh^2(n)). \tag{2.18}$$

We substitute the estimate (2.18) into Theorem 1 [4] and obtain:

$$M[f(i)f(j)] - M[f(i)]M[f(j)] = \frac{o(nh^2(n))}{n^2(n-1)} + O(1/n) = O(1/n). \qquad (2.19)$$

It follows from (2.19) that the estimate of Theorem 1 does not change; therefore the estimate (2.16) of Theorem 2 also does not change.

Taking into account (2.16), we can formulate the following equivalent formulation of the Riemann hypothesis –the boundary of the nontrivial zeros of Riemann zeta function is equal to the order of the standard deviation of Mertens and Liouville arithmetic functions.

Now let us estimate the average value of the above arithmetic functions. First we prove the theorem.

Theorem 3

The summation arithmetic functions of Mertens and Liouville represent an arbitrary symmetric walk at $n \to \infty$. The summative arithmetic function of Liouville is a simple symmetric walk at $n \to \infty$.

Proof

We call an arbitrary symmetric walk the summation arithmetic function $S(n) = \sum_{k=1}^{n} f(k)$, where $f(k)$ is arbitrary independent arithmetic function with zero average value.

The summation functions of Mertens and Liouville have the form $S(n) = \sum_{k=1}^{n} f(k)$. Based on (4) the arithmetic functions $f(k)$ are independent for the value $n \to \infty$ and based on Corollaries 1 and 2 of Theorem 1 [2] and Theorems 1 have zero average value. Therefore, Mertens and Liouville functions represent an arbitrary symmetric walk.

We shall consider a simple symmetric walk the summation arithmetic function $S(n) = \sum_{k=1}^{n} f(k)$, where $f(k)$ are independent arithmetic functions taking only two values $f(k) = 1$ and $f(k) = -1$ which probability is equal to $1/2$.

The summation function of Liouville has the form $S(n) = \sum_{k=1}^{n} f(k)$. Based on (4), the arithmetic functions $f(k)$ are independent for the value $n \to \infty$ and based on Corollaries 1 and 2 of Theorem 1 [2] and Theorem 1 take only two values $f(k)=1$ and $f(k)=-1$, which probability is equal to $1/2$. Therefore, Liouville summation function is a simple walk.

However, despite the fact that the average values of Mobius and Liouville arithmetic functions are zero for the value $n \to \infty$, the average value of Mertens and Liouville summation function are nonzero, and, almost everywhere, they are an unbounded function. Let us prove this for Liouville summation function.

Theorem 4

The average value of the summation arithmetic function of Liouville $M[L(n)]$, almost everywhere, is an unbounded function.

Proof

$$M[L(n)] = \frac{L(1)+L(2)+...+L(n)}{n} = \frac{nf(1)}{n} + \frac{(n-1)f(2)}{n} + ... + \frac{f(n)}{n}.$$

Based on Theorem 3, the summation arithmetic function of Liouville $L(n) = \sum_{k=1}^{n} f(k)$ is a simple symmetric walk.

The following assertion is true (page 373 Example [5]) - series $\sum f(k)a_k$ (where $|a_k| \leq c$) converges, almost everywhere, if and only if $\sum a_k^2 < \infty$ and $f(k)$ is a simple symmetric walk.

The condition for the boundedness of the terms of the series is satisfied in our case - $|a_k| = |\frac{n-k}{n}| \leq 1$.

On other hand:

$$\sum_{k=0}^{n-1}(1-\frac{k}{n})^2 = \frac{1^2+2^2+...+n^2}{n^2} = \frac{n(n+1)(2n+1)}{6n^2}.$$

This sequence tends to $\infty$. when the value $n \to \infty$ Therefore, the series $\sum a_k^2$ diverges and the series $\sum f(k)a_k$ diverges almost everywhere.

Consequently, the average value of the summation arithmetic function of Liouville $M[L(n)]$ almost everywhere is an unbounded function.

The question arises - what is the order of growth of the average value of the above summation arithmetic functions.

Theorem 5

The following estimate is performed for the average value of Liouville summation function:

$$M[L(n)] \leq Cn^{1/2} \log(n) \text{ or } M[L(n)] = O(n^{1/2} \log(n)), \qquad (2.20)$$

where $C$ is a constant.

Proof

It is known that the value $L(n) < 0$ in the interval $1 < n < 906150257$.

However, Liouville summation function has an infinite number of zeros at $n \geq 906150257$.

It is even proved in [6] that:

$$\lim_{n \to \infty} \sup(L(n)/\sqrt{n}) > 0,0618672. \qquad (2.21)$$

It is proved [4] that the standard deviation of Liouville summation arithmetic function $L(n)$ has the order $\sigma[L(n)] = O(n^{1/2})$. Therefore, if we assume that $M[L(n)] > Cn^{1/2} \log(n)$, then starting from a certain value $m$, for all $n > m$, the function $L(n)$ would not have more zeros, but would shift to negative values, which contradicts [6]. Consequently, this assumption is false and the estimate (2.20) is true.

The estimate (2.20) corresponds to [7], where the following approximation of the average value of Liouville summation function is given:

$$L_0(n) = 1 - 0,684765\sqrt{n}.$$

Approximation of the average value of Mertens arithmetic function also has the form:

$$M_0(n) \leq 1 - C\sqrt{n},$$

where $C < 0,0684765$.

## 3. ASYMPTOTIC BEHAVIOR OF MERTENS AND LIOUVILLE SUMMATION ARITHMETIC FUNCTIONS

Theorem 6

The order of Mertens and Liouville summation arithmetic functions $S(n)$ almost everywhere is equal to:

$$S(n) = O(M[S(n)]) + O(n^{1/2} \varphi(n)), \qquad (3.1)$$

where $\varphi(n)$ is any slowly increasing function.

Proof

It is known [1] that almost everywhere is satisfied for an arithmetic function $S(n)$:

$$|S(n) - M[S(n)]| \leq \varphi(n)\sigma[S(n)],$$

where $\varphi(n)$ is any slowly increasing function.

Based on [2.16], this inequality can be written for Mertens and Liouville summation arithmetic functions in the form:

$$|S(n) - M[S(n)]| \leq \varphi(n) n^{1/2}. \qquad (3.2)$$

Having in mind (3.2), the order, almost everywhere, is equal to:

$$S(n) = O(M[S(n)]) + O(n^{1/2} \varphi(n)),$$

which corresponds to (3.1).

Consequence

Based on Theorems 5 and 6, for Liouville summation arithmetic function, almost everywhere, the following estimate holds:

$$L(n) = O(n^{1/2} \log(n)) + O(n^{1/2}\varphi(n)) = O(n^{1/2} \log(n)), \qquad (3.3)$$

which is equivalent to Riemann hypothesis.

It was shown in [8] that the following estimate holds for Mertens function:

$$\lim_{n\to\infty}\sup(\mathrm{M}(n)/\sqrt{n}) > 1,06,$$

which is higher than (2.21). Thus, it is possible that the order of growth of Mertens function is higher than for Liouville summations function.

Therefore, the question arises as to the difference in the order of growth of Mertens and Liouville summation arithmetic functions. To answer this question, we use the methods of complex integration.

Perron's formula [9] for summation arithmetic functions has the form:

$$\sum_{n\leq x} f(n) = \frac{1}{2\pi i}\int_{b-iT}^{b+iT} F(s)\frac{x^s}{s}ds + R_1(x), \qquad (3.4)$$

where $b > a, x \geq 2, T \geq 2, F(s) = \sum_{n=1}^{\infty}\frac{f(n)}{n^s}, \alpha > 0$ are determined from the condition $\sum_{n=1}^{\infty}\frac{|f(n)|}{n^\sigma} \ll \frac{1}{(\sigma-a)^\alpha}$, where $\sigma \to a+0$, and for the function $R_1(x)$ the following estimate holds:

$$R_1(x) \ll \frac{x^b}{T(b-a)^\alpha} + 2^b\left(\frac{x\log x}{T} + \log\frac{T}{b} + 1\right)\max_{x/2\leq n\leq 3x/2}|f(n)|. \qquad (3.5)$$

We consider a special case of Perron's formula for summation arithmetic functions $M(x), L(x)$.

We take in (3.5) the values; $a = 1, b = 1 + 1/\log x$ and consider the integral along the contour $\Gamma$, which is a rectangle; $b \pm iT, \beta \pm iT$, where $\beta$ is the boundary of the nontrivial zeros of the Riemann zeta-function $\zeta(s)$.

This is due to the fact that the data of the summation arithmetic functions $\zeta(s)$ is in the denominator $F(s)$ of formula (3.4). Therefore, we choose the value $\alpha = 1$ and the integration contour, which does not capture the non-trivial zeros of the function $\zeta(s)$.

Now we consider the integral over the given contour $\Gamma$:

$$\frac{1}{2\pi i}\int_{\Gamma} F(s)\frac{x^s}{s}ds = \frac{1}{2\pi i}\int_{b-it}^{b+it} F(s)\frac{x^s}{s}ds + I_1 + I_2 + I_3, \tag{3.6}$$

where $I_1 = \dfrac{1}{2\pi i}\int_{\beta+it}^{b+it} F(s)\dfrac{x^s}{s}ds = \dfrac{1}{2\pi i}\int_{\Gamma_1} F(s)\dfrac{x^s}{s}ds$, $I_1 = \dfrac{1}{2\pi i}\int_{\beta+it}^{b+it} F(s)\dfrac{x^s}{s}ds = \dfrac{1}{2\pi i}\int_{\Gamma_1} F(s)\dfrac{x^s}{s}ds$

$$I_3 = \frac{1}{2\pi i}\int_{\beta-it}^{\beta+it} F(s)\frac{x^s}{s}ds = \frac{1}{2\pi i}\int_{\Gamma_3} F(s)\frac{x^s}{s}ds.$$

It follows from (3.6):

$$\frac{1}{2\pi i}\int_{b-it}^{b+it} F(s)\frac{x^s}{s}ds = \frac{1}{2\pi i}\int_{\Gamma} F(s)\frac{x^s}{s}ds + I_1 + I_2 + I_3. \tag{3.7}$$

We substitute (3.7) into (3.4):

$$\sum_{n\leq x} f(n) = \frac{1}{2\pi i}\int_{\Gamma} F(s)\frac{x^s}{s}ds + I_1 + I_2 + I_3 + R_1(x). \tag{3.8}$$

Denote by $R_2(x) = I_1 + I_2 + I_3$ and having in mind (3.8) we obtain:

$$\sum_{n\leq x} f(n) = \frac{1}{2\pi i}\int_{\Gamma} F(s)\frac{x^s}{s}ds + R(x), \tag{3.9}$$

where $R(x) = R_1(x) + R_2(x)$.

The following relation holds for Mertens and Liouville summations arithmetic functions:

$$\sum_{n\leq x} f(n) = o(x). \tag{3.10}$$

Thus, if this condition is fulfilled, then having in mind (3.9), (3.10) we obtain:

$$\frac{1}{2\pi i}\int_{\Gamma} F(s)\frac{x^s}{s}ds = 0, \tag{3.11}$$

$$R(x) = R_1(x) + R_2(x) = o(x). \tag{3.12}$$

Consequently, the summation arithmetic functions: $M(x), L(x)$ satisfy formulas (3.11), (3.12).

Indeed, for example, for $M(x)$ the function $F(s) = \dfrac{1}{\zeta(s)}$. Since $\zeta(s) \neq 0$ inside $\Gamma$, then $F(s)\dfrac{x^s}{s}$ is regular inside $\Gamma$ and (3.11) is satisfied.

Theorem 7

The order of growth of the ratio of the summed arithmetic functions of Mertens and Liouville is:

$$M(n)/L(n) = O(\log^A n),$$

where $A > 0$.

Proof

It is performed $\max_{x/2 \leq n \leq 3x/2} |f(n)| \leq 1$ for the summation arithmetic functions $M(x), L(\mathrm{x})$, therefore, based on (3.5), we obtain:

$$R_1(x) = O(\dfrac{x \log x}{T}). \qquad (3.13)$$

It is true $F(s) = 1/\zeta(s)$ for Mertens function. Therefore, we take the contour of integration for Mertens function $M(\mathrm{x})$, as in the particular case of Perron's formula considered above.

Based on [10], the following estimate is performed on the horizontal sections $\Gamma 1$ and $\Gamma 2$:

$$\int_{\Gamma_1, \Gamma_2} \dfrac{x^s}{s\zeta(s)} ds = O(\dfrac{\log^B T}{T} \int_\beta^b x^\sigma d\sigma) = O(\dfrac{x \log^B T}{T \log x}). \qquad (3.14)$$

The following estimate is carried out on the vertical section $\Gamma 3$:

$$\int_{\beta - it}^{\beta + it} \dfrac{x^s}{s\zeta(s)} ds = O(\log^B T \cdot x^\beta \cdot \int_2^T \dfrac{dt}{t}) = O(\log^{B+1} T \cdot x^\beta). \qquad (3.15)$$

Here we use the boundary of non-trivial zeros $\beta$ instead of the value $1 - A/\log T$ of the Vallee-Poussin theorem on page 83 [10]

Having in mind (3.14) and (3.15) the following estimate holds for Mertens function:

$$R_2(x) = O(\frac{x \log^B T}{T \log x} + \log^{B+1} T \cdot x^\beta). \qquad (3.16)$$

Consequently, based on (3.13) and (3.16), the following estimate for Mertens function holds:

$$R(x) = R_1(x) + R_2(x) = O(\frac{x \log x}{T} + \frac{x \log^B T}{T \log x} + \log^{B+1} T \cdot x^\beta). \qquad (3.17)$$

We put the value $T = x^{1/2}$ in (3.17), then for any $x \geq 2$, we get:

$$M(x) = \sum_{n \leq x} \mu(n) = R(x) = O(x^{1/2} \log x + x^{1/2} \log^{B-1} x + x^\beta \log^{B+1} x), \qquad (3.18)$$

where $B > 0$.

Based on (3.18), taking into account that $\beta \geq 1/2$ we obtain:

$$M(x) = R(x) = O(x^\beta \log^{B+1} x). \qquad (3.19)$$

Now we consider Liouville summation function.

It is performed $F(s) = \zeta(2s)/\zeta(s)$ for Liouville summation function.

We take the integration contour for the function $L(x)$ as in the particular case of Perron's formula considered above.

Based on properties of the function $F(s) = \zeta(2s)/\zeta(s)$, we similarly obtain:

$$R_2(x) = \int_{\Gamma_1, \Gamma_2, \Gamma_3} \frac{\zeta(2s) x^s}{s \zeta(s)} ds = O(x^\beta \cdot \log T). \qquad (3.20)$$

Consequently, having in mind (3.13) and (3.20), the following estimate for the function $L(x)$ holds:

$$R(x) = R_1(x) + R_2(x) = O(\frac{x \log x}{T} + x^\beta \cdot \log T). \qquad (3.21)$$

Suppose that $T = x^{1/2}$, then based on (3.21), for any $x \geq 2$, we get:

$$L(x) = \sum_{n \leq x} \lambda(n) = O(x^{1/2} \log(x) + x^{\beta} \bullet \log(x)). \tag{3.22}$$

Based on (3.22), taking into account that $\beta \geq 1/2$ we obtain:

$$L(x) = R(x) = O(x^{\beta} \log(x)). \tag{3.23}$$

We define the order of growth of the ratio of these functions having in mind (3.19) and (3.23):

$$M(x)/L(x) = O(\log^{B}(x)). \tag{3.24}$$

Consequence

Having in mind (3.3) and (3.24), we obtain that, almost everywhere, the following estimate for Mertens summation arithmetic function is valid:

$$M(n) = O(n^{1/2} \log^{B+1} n), \tag{3.25}$$

which is equivalent to Riemann hypothesis.

## 4. CONCLUSION AND SUGGESTIONS FOR FURTHER WORK

The next article will continue to study the behavior of arithmetic functions.

## 5. ACKNOWLEDGEMENTS

Thanks to everyone who has contributed to the discussion of this paper. I am grateful to everyone who expressed their suggestions and comments in the course of this work.